\documentclass[12pt]{article}
\usepackage{ucs}
\usepackage[utf8x]{inputenc}
\usepackage{amsfonts,amssymb}
\usepackage{amsmath}
\usepackage{amsthm}
\usepackage[english]{babel}
\usepackage{indentfirst}
\usepackage{geometry}
\usepackage{graphicx}
\usepackage{multicol}

\newtheorem{Th}{Theorem}
\newtheorem{Cor}{Corollary}

\title{Diophantine approximations with Pisot numbers}
\author{V. Zhuravleva\thanks{The research was supported by RNF grant 14-11-00433}}
\date{}

\begin{document}
\maketitle

\begin{abstract}
Let $\alpha$ be a Pisot number. Let $L(\alpha)$ be the largest positive number such that for some $\xi=\xi(\alpha)\in \mathbb R$ the limit points of the sequence of fractional parts $\{\xi \alpha^n\}_{n=1}^{\infty}$ all lie in the interval $[L(\alpha), 1-L(\alpha)]$. In this paper we show that if $\alpha$ is of degree at most 4 or $\alpha\le \frac{\sqrt 5 + 1}{2}$ then $L(\alpha)\ge \frac{3}{17}$. Also we find explicitly the value of $L(\alpha)$ for certain Pisot numbers of degree 3.
\end{abstract}

\section{Introduction}

Our introduction is organised as follows. In subsections 1.1-1.4 below we discuss the definitions and known results, related to Diophantine approximations with Pisot numbers. Our new results are formulated in subsection 1.5.

\subsection{Pisot numbers}

A real algebraic integer $\alpha$ is called a Pisot number if it is greater than 1 and all its remaining conjugates lie strictly inside the unit circle $\mathbb D = \{z \in \mathbb C : |z|<1\}$. Pisot numbers appeared in A. Thue's paper in 1912 \cite{ref1}. Basic properties of these numbers were independetly studied by Ch. Pisot \cite{ref2} and T. Vijayaraghavan \cite{ref3}. In J. W. S. Cassels's book \cite{ref4} these numbers are mentioned as Pisot–Vijayaraghavan numbers.

By
$$
P(x)=x^d - a_1 x^{d-1} - a_2 x^{d-2} - \ldots - a_k x^{d-k} - \ldots a_{d-1} x - a_{d}, \quad \qquad \text {where $a_i \in \mathbb Z$} 
$$
we denote the minimal polynomial of $\alpha$.

There is no rigorous description of the set of coefficients $(a_1, \ldots, a_d)$ that correspond to Pisot numbers. Rouche's theorem leads to the following sufficient condition. 

\noindent {\bf Proposition 1.1} [\cite{ref5}, ch. 5].
{\it If $1+\sum_{i=1}^d a_i<0$ and $|a_1|>1+\sum_{i=2}^{d}|a_i|$, then the maximal root of $P(X)$ is a Pisot number.}

However the most famous Pisot number $\alpha=\frac{\sqrt 5 +1}{2}$ doesn't satisfy these inequalities.

The set of Pisot numbers less than $\frac{\sqrt 5 + 1}{2}$ is described completely. We have the following statement.

\noindent {\bf Proposition 1.2} [\cite{ref5}, ch. 7].
{\it Pisot numbers less than $\frac{\sqrt 5+1}{2}$ can be arranged in increasing order as
$$\theta_2=\theta'_1<\theta'_2<\theta'_3<\theta_3<\theta'_4<\theta_4<\theta'_5<\theta''<\theta_5<\theta'_6<\ldots<$$ $$<\theta'_k<\theta_k<\ldots<\frac{\sqrt 5 +1}{2}.$$
The algebraic integers $\theta_{2p}$ (respectively $\theta_{2p+1}$) are roots of the polynomials 
$$T_{2p}=\frac{x^{2p}(x^2-x-1)+1}{x-1}, \quad (\text{{\it respectively} } T_{2p+1}=\frac{x^{2p+1}(x^2-x-1)+1}{x^2-1}).$$
The Pisot numbers $\theta''$, $\theta'_k$ are zeros respectively of the polynomials
$$T''=1-x+x^2-x^4+2x^5-x^6, \quad T'_k=1-x^2+x^k(1+x-x^2).$$}

Two smallest Pisot numbers ($\theta_1$, $\theta_2$ from Propostion 1.2) were found by C. L. Siegel \cite{ref6}.

\subsection{Linear recurrence sequences}

Consider a sequence $X=(x_n)_{n=1}^{\infty}$ which satisfies the relation
\begin{equation}
x_{n}=a_1 x_{n-1} + a_2 x_{n-2} + \ldots + a_k x_{n-k} + \ldots + a_{d-1} x_{n-d+1} + a_{d} x_{n-d}. 
\end{equation}

In the present paper we study the value of
$$
L(X):=\sup_{\xi \in \mathbb R} \liminf_{n \to \infty} ||\xi x_n||
$$
(here $||\cdot||$ denotes the distance to the nearest integer, while $\{\cdot\}$ stands for the fractional part, and  $[\cdot]$ denotes the integral part).

It is obvious that if two sequences $X=(x_n)_{n=1}^{\infty}$, $Y=(y_n)_{n=1}^{\infty}$ satisfy (1) and $x_n \nrightarrow 0$, $y_n \nrightarrow 0$ as $n \rightarrow \infty$, then
$$L(X)=L(Y).$$

So, if $X=(x_n)_{n=1}^{\infty}$ satisfies (1), then
$$
L(X)=\begin{cases}
0, \text{\quad\quad\quad\quad\quad if $x_n \rightarrow 0$ as $n \rightarrow \infty$,}\\
L(\alpha)>0, \text{ \quad if $x_n \nrightarrow 0$ as $n \rightarrow \infty$.}
\end{cases}
$$

The inequality $L(\alpha)>0$ holds since the sequence $X$ is lacunary when $x_n \nrightarrow 0$ as $n \rightarrow \infty$. We call a sequence $X$ lacunary if there exists $\lambda > 1$ such that for all sufficiently large $n$ we have $|x_{n+1}|>\lambda|x_n|.$ For any lacunary secuence $X$ A.Y. Khintchin \cite{ref7} proved that there exist $\xi \in \mathbb R$ and $\gamma>0$ such that for all $n \in \mathbb N$ we have
$$
||\xi x_n||\ge \gamma.
$$

This theorem has an interesting history (see \cite{Peres}, \cite{ref8}, \cite{ref9}, \cite{ref10}, \cite{Rochev}).

\subsection{Value of $L(\alpha)$ for Pisot numbers}

In some cases the value of $L(\alpha)$ can be calculated explicitely. Here we list all known results. 

\noindent{\bf Theorem A.} {\it We have the following equalities.

1) If $\alpha=\frac{\sqrt 5 +1}{2}$, then $L(\alpha)=\frac{1}{5}$;

2) If $\alpha$ is the maximal zero of $x^3=x+1$, then $L(\alpha)=\frac{1}{5}$;

3) If $\alpha$ is the maximal zero of $x^4=x^3+1$, then $L(\alpha)=\frac{3}{17}$.}

In fact, this theorem represents the major results of two author's papers \cite{ref13}, \cite{ref14}. Different proofs were given by A. Dubickas in \cite{ref12} and \cite{newdub}.

The following statement appeared in Dubickas' paper \cite{ref12}.

\noindent{\bf Theorem B.} {\it Let $\alpha$ be a Pisot number and let $P(X)$ be its minimal polynomial. Suppose that $\sum_{i=1}^{d} a_i$ is odd. Then $L(\alpha)=\frac{1}{2}.$}

The following two author's theorems from \cite{ref15} generalize Theorem 1.2 and 1.3 from Dubickas' paper \cite{ref12}.

\noindent {\bf Theorem C.} {\it Let $\alpha$ be a Pisot number and let $P(X)$ be its minimal polynomial. Suppose that $a_i\ge0$, if $i$ is odd, $a_i\le0$, if  $i$ is even and $\sum_{i=1}^{d} a_i$ is even. Then $$L(\alpha)=\frac {\sum_{i=1}^{d} {|a_i|}}{2\sum_{i=1}^{d} {|a_i|}+2}.$$}

\noindent {\bf Theorem D.} {\it Let $\alpha$ be a Pisot number and let $P(X)$ be its minimal polynomial. Suppose that $a_i \in \mathbb N$, $i=1, \ldots, d$ and $\sum_{i=1}^{d} a_i$ is even. Also let $d>1$ and $\alpha\ne\frac{\sqrt 5 + 1}{2}$. Then $$L(\alpha)=\frac{\sum_{i=1}^{d} a_i-2}{2\sum_{i=1}^{d} a_i-2}.$$}

(In \cite{ref12} Dubickas has these results for Pisot numbers of degree $\le 2$).

\subsection{Reduced length of a polynomial}

In \cite{ref11} Dubickas introduced the reduced length $d(P)$ of a polynomial. Consider the set $\Gamma=\{g_0+g_1 x+\ldots+g_m x^m \in \mathbb R[x]$, where $g_0=1$ and $g_m=1\}$. The reduced length of $P$ is defined as

$$
d(P):=\inf _{G \in \Gamma} D(PG),
$$
where $D(\cdot)$ is usual length of a polynomial, which is the sum of the absolute values of its coefficients. Then we have the following statement. 

\noindent {\bf Theorem E \cite{ref11}}. {\it Let $\alpha$ be a Pisot number. If $\xi \notin \mathbb Q(\alpha)$, then $$\limsup_{n \to \infty} \{\xi \alpha^n \}-\liminf_{n \to \infty}\{\xi \alpha^n \}\ge \frac{1}{d(P)},$$ where $P$ is the minimal polynomial for $\alpha$.}

This theorem immediately leads to the following

\noindent {\bf Corollary.} {\it If $\xi \notin \mathbb Q(\alpha)$, then }$$\liminf_{n \to \infty} ||\xi x_n||\le\frac{1}{2}-\frac{1}{2d(P)}\le \frac{1}{2}-\frac{1}{2D(P)}=\frac {\sum_{i=1}^{d} {|a_i|}}{2\sum_{i=1}^{d} {|a_i|}+2}.$$

\subsection{New results}

In this section we formulate our new results. Theorem 1 below gives the simplest upper bound for  $L(\alpha)$, which is similar to Dubickas' one in the case $d=1$.

\begin{Th} 
Let $\alpha$ be a Pisot number and let $P(X)$ be its minimal polynomial. Suppose that $\sum_{i=1}^{d} a_i$ is even. Then $$L(\alpha)\le\frac {\sum_{i=1}^{d} {|a_i|}}{2\sum_{i=1}^{d} {|a_i|}+2}.$$
\end{Th}

Note that from Theorem C we know that this upper bound can't be improved without adding supplementary conditions.

Put
$$
s_{j,t}=\sum_{\substack {i=1  \\ i \equiv j\pmod t }}^{d}a_i, \qquad \text {if } j=1 \ldots t.
$$

In particular, $$s_{1,1}=\sum_{i=1}^{d}a_i,$$ $$s_{1,2}=\sum_{\substack{i=1 \\ i \text{ is odd}}}^{d}a_i, \qquad s_{2,2}=\sum_{\substack{i=1 \\ i \text{ is even}}}^{d}a_i.$$

The following statement gives us a general lower bound for $L(\alpha)$.

\begin{Th}
Let $\alpha$ be a Pisot number and let $P(X)$ be its minimal polynomial. Suppose that $\sum_{i=1}^{d} a_i$ is even. Then
\begin{multline}
L(\alpha)\ge \max \{ \left \|\frac{s_{1,1}-2}{2s_{1,1}-2}\right \|, \left \| \frac{s_{1,2}-s_{2,2}}{2(s_{1,2}-s_{2,2})+2} \right \|, 
\\ \frac {(s_{1,4}-s_{3,4})^2-|s_{1,4}-s_{3,4}|+(s_{2,4}-s_{4,4}+1)^2-|s_{2,4}-s_{4,4}+1|}{2(s_{1,4}-s_{3,4})^2+2(s_{2,4}-s_{4,4}+1)^2} \}.
\end{multline}
\end{Th}

To formulate our next general result we need a special result concerning just one particular number. Recall that we obtained the results of parts 2, 3 of Theorem A by computer calculations \cite{ref14}. Quite similar calculations lead to Theorem 3, which we would like to anounce here without a complete proof. Some ideas related to our calculations and further computational results are discussed in Sections 7.1, 7.2 below.

\begin{Th}
Let $\alpha$ be the root of $x^3-2x^2-x+1=0$. Then $L(\alpha)= \frac{4}{13}$.
\end{Th}

The statement of Theorem 3 can be obtained by numerical calculations completely described in \cite{ref14}.

Our main result concerning Pisot numbers of degree not greater than 4 is as follows.

\begin{Th}
If $\alpha$ is a Pisot number of degree not greater than 4 and $\alpha$ is not the root of $x^2-x-1=0$, $x^3-x-1=0$, $x^4-x^3-1=0$ or $x^3-2x^2-x+1=0$ then $L(\alpha)\ge \frac13$.
\end{Th}

Theorems 3,4 and A lead to the following corollaries.

\begin{Cor} We have the following statements.

1) If $\alpha$ is a Pisot number of degree $\le 3$, then
$L(\alpha)\ge \frac{1}{5};$

2) If $\alpha$ is a Pisot number of degree $\le 4$, then
$L(\alpha)\ge \frac{3}{17}.$
\end{Cor}

\begin{Cor}
Let $\alpha$ be a Pisot number of degree not greater than 4. Then

1) $L(\alpha)=\frac{3}{17}$ if and only if $\alpha$ is the root of $x^4-x^3-1=0$,

2) $L(\alpha)=\frac{1}{5}$ if and only if $\alpha$ is the root of $x^2-x-1=0$ or $x^3-x-1=0$,

3) $L(\alpha)=\frac{4}{13}$ if and only if $\alpha$ is the root of $x^3-2x^2-x+1=0$.
\end{Cor}

Proposition 1.2 from subsection 1.1 describes all Pisot numbers less than $\frac{\sqrt 5 +1}{2}$. Our approach gives the following lower bound for $L(\alpha)$ for all these numbers.

\begin{Th}
Let $\alpha$ be a Pisot number less than $\frac{\sqrt 5 + 1}{2}$. Then $$L(\alpha)\ge \frac{3}{17}.$$
\end{Th}

Theorem A shows that the upper bounds of Theorems 3 and 4 are optimal.

Using Theorem 2 and Dubickas' methods from \cite{ref12} we obtain exact values of $L(\alpha)$ for certain Pisot numbers of degree 3. These values are presented below.

\begin{Th} Let $\alpha$ be a Pisot number with the minimal polynomial $P(X)=x^3-a_1x^2-a_2x-a_3$. Then we have the following statements.

1) If $a_2<0, a_3<0$, and $a_2+\frac {a_3^2-a_3(a_1-1)+1}{a_1+1}\le 0$

then $$L(\alpha) = \frac {a_1-a_2+a_3}{2(a_1-a_2+a_3+1)};$$ 

2) If $a_2> 0, a_3<0$, and $a_2-\frac {a_3^2-a_3(a_1+1)+1}{a_1-1}\ge 0$
then
$$ L(\alpha) = \frac {a_1+a_2+a_3-2}{2(a_1+a_2+a_3-1)}.$$
\end{Th}

One can see that Theorems C, D and 6 in some cases give the exact values of $L(\alpha)$ for Pisot numbers of degree 3. For some Pisot numbers we can find the exact values $L(\alpha)$ by using computer calculations. In section 7 we formulate some results of such calculations for some Pisot numbers of degree 3 and 4.

Our paper is organised as follows. In Section 2 we give a proof of Theorem 1. Section 3 is devoted to periodic sequences modulo 1, which help us to obtain the lower bounds for $L(\alpha)$ from Theorem 2. We apply these bounds in Sections 4 and 5, where we prove Theorems 4 and 5. Section 6 is devoted to detailed study of Pisot numbers of degree 3. In fact, we find the set of all coefficients $(a_1, a_2, a_3)$ which correspond to Pisot numbers of degree 3. As mentioned above for some of them the exact value of $L(\alpha)$ can be found. In Section 7 we talk about computer methods.

\section {Proof of Theorem 1}

Theorem 1 is proved similarly to Theorem 1.2 from \cite{ref12}. But for the sake of completeness we present its detailed proof.

Put $u_n:=[\xi x_n]$ and $v_n:=\{\xi x_n\}$. Then $x_n=u_n+v_n$. We suppose that there exists $n_0 \in \mathbb N$ such that for all $n>n_0$
$$
v_n \in [k, 1-k], \qquad \text{where } k>\frac {\sum_{i=1}^{d} {|a_i|}}{2\sum_{i=1}^{d} {|a_i|}+2}.
$$

It's clear that for all natural $n>d$ we have
$$
u_n+v_n=\sum_{i=1}^{d}a_i u_{n-i}+\sum_{i=1}^{d}a_i v_{n-i}.
$$

Put
$$
\delta_n:=\sum_{i=1}^{d}a_i v_{n-i}-v_n=u_n-\sum_{i=1}^{d}a_i u_{n-i} \in \mathbb Z.
$$

Let
$$
I_{+}= \left\{i : a_i \ge 0 \right\} \qquad \text{and} \qquad I_{-}= \left\{i : a_i < 0 \right\}.
$$

Therefore
$$
\delta_n\in \left(\sum_{i \in I_{-}} a_i - 1, \sum_{i \in I_{+}} a_i \right).
$$

So 
$$
\delta_n\in \left \{\sum_{i \in I_{-}} a_i , \ldots,  \sum_{i \in I_{+}} a_i -1 \right \},
$$

and one of the following conditions holds.

1) For infinitely many $n\in \mathbb N$ we have
$$
\delta_n \le \frac {\sum_{i=1}^d a_i}{2}-1.
$$

2) For infinitely many $n\in \mathbb N$ we have
$$
\delta_n \ge \frac {\sum_{i=1}^d a_i}{2}.
$$

We consider the first case. Then for infinitely many $n\in \mathbb N$ we have
$$
\sum_{i=1}^{d} a_i v_{n-i}\ge k \sum_{i \in I_{+}} a_i + (1-k)\sum_{i \in I_{-}} a_i=k\sum_{i=1}^d |a_i|+\sum_{i \in I_{-}} a_i.
$$

On the other hand,
$$
\sum_{i=1}^{d} a_i v_{n-i}=v_n+\delta_n\le 1-k + \frac {\sum_{i=1}^d a_i}{2}-1=\frac {\sum_{i=1}^d a_i}{2}-k.
$$

We obtain that
$$
k\sum_{i=1}^d |a_i|+\sum_{i \in I_{-}} a_i \le \frac {\sum_{i=1}^d a_i}{2}-k.
$$

Then
$$
k(\sum_{i=1}^d |a_i|+1) \le \frac {\sum_{i=1}^d |a_i|}{2}.
$$

But we assume that $k>\frac {\sum_{i=1}^{d} {|a_i|}}{2\sum_{i=1}^{d} {|a_i|}+2}$.

We consider the second case. Then for infinitely many $n\in \mathbb N$ we have
$$
\sum_{i=1}^{d} a_i v_{n-i}\le (1-k) \sum_{i \in I_{+}} a_i + k \sum_{i \in I_{-}} a_i=\sum_{i \in I_{+}} a_i - k\sum_{i=1}^d |a_i|.
$$

On the other hand,
$$
\sum_{i=1}^{d} a_i v_{n-i}=v_n+\delta_n\ge k + \frac {\sum_{i=1}^d a_i}{2}.
$$

Then
$$
k+\frac {\sum_{i=1}^d a_i}{2}\le \sum_{i \in I_{+}} a_i - k\sum_{i=1}^d |a_i|.
$$

We obtain the contradiction similar to the one in the first case. This completes the proof.

\section{Periodic sequences modulo 1}

A sequence $X=(x_n)_{n=1}^{\infty}$ is called periodic modulo 1, if there exists an integer $t$ such that  $\{x_n\}=\{x_{n+t}\}$ for all $n \in \mathbb N$. An integer $t$ with such a condition is the length of the period of this sequence. We do not necessarily consider the smallest period. Each sequence $\{x_n\} \subset \mathbb Q$ satisfying (1) has a period. Each period is of the form

\begin{equation}
\frac{r_1}{z}, \frac{r_2}{z}, \ldots, \frac{r_t}{z}, \qquad \text{where } z\in \mathbb N, r_1, \ldots, r_t \in \mathbb Z, 0\le r_1, \ldots, r_t < z.
\end{equation}

Given arbitrary sequence $\{x_n\} \subset \mathbb Q$ and a positive integer $t$ consider $t$ sums
$$
s_{j,t}=\sum_{\substack {i=1  \\ i \equiv j\pmod t }}^{d}a_i, \qquad \text {if } j=1 \ldots t.
$$

The following Lemma comes from the definition of a peiodic sequence.
\\

\noindent {\bf Lemma 3.1 [Criterion of periodic sequence modulo 1].} {\it Consider recurrent relation (1). Let $t$ be positive integer. Then (3) is a $t$-periodic rational sequence satisfying (1) if and only if 
\begin{equation}
\begin{cases}
s_{t,t} r_1+s_{t-1,t}r_2+\ldots + s_{2,t} r_{t-1} +s_{1,t} r_t \equiv r_1 \pmod z \\
s_{t,t} r_2+s_{t-1, t}r_3+\ldots + s_{2,t} r_{t} +s_{1,t} r_1 \equiv r_2 \pmod z \\
\ldots \\
s_{t,t} r_{t-1}+s_{t-1,t} r_t+\ldots + s_{2,t} r_{t-3} +s_{1,t} r_{t-2} \equiv r_{t-1} \pmod z \\
s_{t,t} r_t+s_{t-1, t} r_1+\ldots + s_{2,t} r_{t-2} +s_{1,t} r_{t-1} \equiv r_t \pmod z.
\end{cases}
\end{equation}}

\noindent {\bf Remark.} Conditions (4) depend not on the recurrence (1) itself, but on the sums $s_{i,j}$ only. We take another relation
\begin{equation}
x_{n}=a'_1 x_{n-1} + a'_2 x_{n-2} + \ldots + a'_k x_{n-k} + \ldots + a'_{d'-1} x_{n-d'+1} + a'_{d'} x_{n-d'}. 
\end{equation}

Put $$s'_{j,t}=\sum_{\substack {i=1  \\ i \equiv j\pmod t }}^{d'}a'_i, \qquad \text {if } j=1 \ldots t.$$

Let $s_{j,t}=s'_{j,t}$ for all $j=1, \ldots, t$, and let $z, r_1, \ldots, r_t$ satisfy (4).  We consider periodic modulo 1 sequence $X=(x_n)_{n=1}^{\infty}$ with $\{x_n\}=\frac{r_i}{z}$ if $n \equiv i \pmod t$. Then $X$ satisfy both (1) and (5).

Lemmas 3.2 - 3.4 below can be deduced by easy calculations.

\noindent {\bf Lemma 3.2.} {\it Let $\sum_{i=1}^{d}x_i$ be even. Suppose that $X=(x_n)_{n=1}^{\infty}$ satisfies (1). Also let
$$    
x_i=\frac {s_{1,1}-2}{2s_{1,1}-2},  \quad \text{if $i\in[1,\ldots,d]$.}
$$
Then sequence $X$ is periodic modulo 1. The length of its period is 1.
\\

\noindent {\bf Lemma 3.3.} {\it Let $\sum_{i=1}^{d}x_i$ be even. Suppose that $X=(x_n)_{n=1}^{\infty}$ satisfies (1). Also let
$$    
x_i=\begin{cases}
        \frac {s_{1,2}-s_{2,2}}{2(s_{1,2}-s_{2,2}+1)}, & \text{if $i$ is odd, $i\in[1,\ldots,d]$}; \\
        \frac{s_{1,2}-s_{2,2}+2}{2(s_{1,2}-s_{2,2}+1)}, & \text{if $i$ is even, $i\in[1,\ldots,d]$}.
  \end{cases}
$$
Then sequence $X$ is periodic modulo 1. The length of its period is 2.}
\\

\noindent {\bf Lemma 3.4.} {\it Let $\sum_{i=1}^{d}x_i$ be even. Suppose that  $X=(x_n)_{n=1}^{\infty}$ satisfies (1). Also let
$$
x_i=
\begin{cases}
\frac {(s_{1,4}-s_{3,4})^2-(s_{1,4}-s_{3,4})+(s_{2,4}-s_{4,4}+1)^2-(s_{2,4}-s_{4,4}+1)}{2(s_{1,4}-s_{3,4})^2+2(s_{2,4}-s_{4,4}+1)^2}, & \text{if $i \equiv 1\pmod 4$, $i\in[1,\ldots,d]$;} \\
\frac {(s_{1,4}-s_{3,4})^2-(s_{1,4}-s_{3,4})+(s_{2,4}-s_{4,4}+1)^2+(s_{2,4}-s_{4,4}+1)}{2(s_{1,4}-s_{3,4})^2+2(s_{2,4}-s_{4,4}+1)^2}, & \text{if $i \equiv 2\pmod 4$, $i\in[1,\ldots,d]$;} \\
\frac {(s_{1,4}-s_{3,4})^2+(s_{1,4}-s_{3,4})+(s_{2,4}-s_{4,4}+1)^2+(s_{2,4}-s_{4,4}+1)}{2(s_{1,4}-s_{3,4})^2+2(s_{2,4}-s_{4,4}+1)^2}, & \text{if $i \equiv 3\pmod 4$, $i\in[1,\ldots,d]$;} \\
\frac {(s_{1,4}-s_{3,4})^2+(s_{1,4}-s_{3,4})+(s_{2,4}-s_{4,4}+1)^2-(s_{2,4}-s_{4,4}+1)}{2(s_{1,4}-s_{3,4})^2+2(s_{2,4}-s_{4,4}+1)^2}, & \text{if $i \equiv 4\pmod 4$, $i\in[1,\ldots,d].$}
\end{cases}
$$
Then sequence $X$ is periodic modulo 1. The length of its period is 4.}}
\\

Here we do not want to give the complete proofs of Lemmas 3.2-3.4. All of them are similar. In their proofs we use Lemma 3.1. As an example we prove the first assertion of Lemma 3.4. According to Lemma 3.1. we have to prove 
\begin{multline}
s_{4,4}((s_{1,4}-s_{3,4})^2-(s_{1,4}-s_{3,4})+(s_{2,4}-s_{4,4}+1)^2-(s_{2,4}-s_{4,4}+1)) \\
+s_{3,4}((s_{1,4}-s_{3,4})^2-(s_{1,4}-s_{3,4})+(s_{2,4}-s_{4,4}+1)^2+(s_{2,4}-s_{4,4}+1)) \\
+s_{2,4}((s_{1,4}-s_{3,4})^2+(s_{1,4}-s_{3,4})+(s_{2,4}-s_{4,4}+1)^2+(s_{2,4}-s_{4,4}+1)) \\ 
+s_{1,4}((s_{1,4}-s_{3,4})^2+(s_{1,4}-s_{3,4})+(s_{2,4}-s_{4,4}+1)^2-(s_{2,4}-s_{4,4}+1)) \equiv\\
\equiv(s_{1,4}-s_{3,4})^2-(s_{1,4}-s_{3,4})+(s_{2,4}-s_{4,4}+1)^2-(s_{2,4}-s_{4,4}+1) \pmod {(2(s_{1,4}-s_{3,4})^2+2(s_{2,4}-s_{4,4}+1)^2)}.
\end{multline}

It's easy to see that both hand sides of (6) are even. That's why we need to prove this equality $\mod {(s_{1,4}-s_{3,4})^2+(s_{2,4}-s_{4,4}+1)^2}$. Then (6) is equivalent to
\begin{multline}\notag
(s_{4,4}+s_{3,4}+s_{2,4}+s_{1,4}-1)((s_{1,4}-s_{3,4})^2+(s_{2,4}-s_{4,4}+1)^2) \\
+s_{1,4}^2+s_{2,4}^2+s_{3,4}^2+s_{4,4}^2-2s_{1,4}s_{3,4}-2s_{2,4}s_{4,4}-s_{1,4}+s_{2,4}+s_{3,4}-s_{4,4}\equiv \\
\equiv -s_{1,4}-s_{2,4}+s_{3,4}+s_{4,4}-1 \pmod {((s_{1,4}-s_{3,4})^2+(s_{2,4}-s_{4,4}+1)^2)}.
\end{multline}

Finally, we get $$(s_{4,4}+s_{3,4}+s_{2,4}+s_{1,4})((s_{1,4}-s_{3,4})^2+(s_{2,4}-s_{4,4}+1)^2) \equiv 0 \pmod {((s_{1,4}-s_{3,4})^2+(s_{2,4}-s_{4,4}+1)^2}).$$

So we proved the first case of Lemma 3.4.

In \cite{ref12} Dubickas proved the following statement.
\\

\noindent {\bf Lemma 3.5.} {\it 
Suppose that there exists periodic modulo 1 sequence  $X^*=(x^*_n)_{n=1}^{\infty}$ that satisfies (1). Suppose that the length of the period of this sequence is $t$. Then $$L(\alpha)\ge \min_{i=1, \ldots, t} ||x_i||.$$ }

From Lemmas 3.2-3.5 we get Theorem 2.

\section{Lower bound for Pisot numbers of degree $\le 4$}

The bound from Theorem 2 is useful only if the right hand side of (2) is not $0$.

In the case $d\le 4$ the right hand side of (2) is equal to 0 if and only if 

\begin{equation}
\begin{cases}
a_1+a_2+a_3+a_4=0, \\
a_1+a_3=a_2+a_4, \\
(a_1-a_3)^2-|a_1-a_3|+(a_2-a_4+1)^2-|a_2-a_4+1|=0.
\end{cases}
\end{equation}

or

\begin{equation}
\begin{cases}
a_1+a_2+a_3+a_4=2, \\
a_1+a_3=a_2+a_4, \\
(a_1-a_3)^2-|a_1-a_3|+(a_2-a_4+1)^2-|a_2-a_4+1|=0.
\end{cases}
\end{equation}

or

\begin{equation}
\begin{cases}
a_1+a_2+a_3+a_4=0, \\
a_1+a_3-a_2-a_4=-2, \\
(a_1-a_3)^2-|a_1-a_3|+(a_2-a_4+1)^2-|a_2-a_4+1|=0.
\end{cases}
\end{equation}

or

\begin{equation}
\begin{cases}
a_1+a_2+a_3+a_4=2, \\
a_1+a_3-a_2-a_4=-2, \\
(a_1-a_3)^2-|a_1-a_3|+(a_2-a_4+1)^2-|a_2-a_4+1|=0.
\end{cases}
\end{equation}

In each of systems (7)-(10) the last equation is the same. It is equivalent to $$a_1-a_3 \in \{-1,0,1\}, \quad a_2-a_4+1 \in \{-1,0,1\}.$$

System (7) has 9 solutions $(\frac12,0,-\frac{1}{2},0)$, $(\frac12, -\frac12, -\frac12, \frac12)$, $(\frac12, -1, -\frac12, 1)$, $(0,0,0,0)$, $(0, -\frac12, 0, \frac12)$, $(0,-1,0,1)$, $(-\frac12,0,\frac12,0)$, $(-\frac12,-\frac12,\frac12,\frac12)$, $(-\frac12,-1,\frac12,1)$. 

System (8) has 9 solutions $(1,\frac12,0,\frac12)$, $(1,0,0,1)$, $(1,-\frac12,0,\frac32)$, $(\frac12,\frac12,\frac12,\frac12)$, $(\frac12,0,\frac12,1)$, $(\frac12,-\frac12,\frac12,\frac32)$, $(0,\frac12,1,\frac12)$, $(0,0,1,1)$, $(0,-\frac12,1,\frac32)$.

System (9) has 9 solutions $(0,\frac12,-1,\frac12)$, $(0,0,-1,1)$, $(0,-\frac12,-1,\frac32)$, $(-\frac12,\frac12,-\frac12,\frac12)$, $(-\frac12,0,-\frac12,1)$, $(-\frac12,-\frac12,-\frac12,\frac32)$, $(-1,\frac12,0,\frac12)$, $(-1,0,0,1)$, $(-1,-\frac12,0,\frac32)$.

System (10) has 9 solutions $(\frac12,1,-\frac12,1)$, $(\frac12,\frac12,-\frac12,\frac32)$, $(\frac12,0,-\frac12,2)$, $(0,1,0,1)$, $(0,\frac12,0,\frac32)$, $(0,0,0,2)$, $(-\frac12,1,\frac12,1)$, $(-\frac12,\frac12,\frac12,\frac32)$, $(-\frac12,0,\frac12,2)$.

Among these solutions only $(1,0,0,1)$ correspond to a Pisot number. Namely, it is the maximal root of equation $x^4=x^3+1$. For this root we know that $L(\alpha)=\frac{3}{17}$ (\cite{ref14}). 

Now we prove that if $\alpha$ is a Pisot number of degree not greater than 4 and $\alpha$ is not the root of $x^4-x^3-1=0$ then

\begin{multline}
L(\alpha)\ge \max \{ \left \| \frac{a_1+a_2+a_3+a_4-2}{2(a_1+a_2+a_3+a_4-1)} \right \|, \left \| \frac{a_1-a_2+a_3-a_4}{2(a_1-a_2+a_3-a_4+1)} \right \|,  
\\ \frac {(a_1-a_3)^2-|a_1-a_3|+(a_2-a_4+1)^2-|a_2-a_4+1|}{2(a_1-a_3)^2+2(a_2-a_4+1)^2}\}\ge \frac{1}{5}. 
\end{multline}

All Pisot numbers of degree not greater than 4 except the zero of $x^4-x^3-1=0$ do not satisfy any of systems (7)-(10).

If $a_1+a_2+a_3+a_4$ is odd then $L(\alpha)=\frac{1}{2}$. So we consider only the cases when $a_1+a_2+a_3+a_4$ is even. 

If 
\begin{equation}
a_1+a_2+a_3+a_4\ge 4 \text{ or } a_1+a_2+a_3+a_4\le-2, 
\end{equation}
then $$\frac{a_1+a_2+a_3+a_4}{2(a_1+a_2+a_3+a_4-1)}\ge \frac{1}{3}.$$

If 
\begin{equation}
a_1-a_2+a_3-a_4\ge 2 \text{ or } a_1-a_2+a_3-a_4\le -4\, 
\end{equation}
then $$ \frac{a_1-a_2+a_3-a_4}{2(a_1-a_2+a_3-a_4+1)} \ge \frac{1}{3}.$$

If both of (12) and (13) are not satisfied we have 

\begin{equation}
a_1+a_2+a_3+a_4\in\{0,2\} \text{ and } a_1-a_2+a_3-a_4\in\{-2,0\}.
\end{equation}

Then we must consider the third bound from (11). 

Put $A:=|a_1-a_3|$ and $B:=|a_2-a_4+1|$. We note that $A$ and $B$ have different parities. We have (14) and we don't consider solutions of (7)-(10). This gives $(A,B)\ne (0,1), (1,0)$. So $(A,B)$ with the smallest possible value of $A^2+B^2$ are $(2,1)$ and $(1,2)$. Then it's easy to see that 

\begin{multline}
z:=\frac {(a_1-a_3)^2-|a_1-a_3|+(a_2-a_4+1)^2-|a_2-a_4+1|}{2(a_1-a_3)^2+2(a_2-a_4+1)^2}=\\
=\frac{1}{2}-\frac{|a_1-a_3|+|a_2-a_4+1|}{2(a_1-a_3)^2+2(a_2-a_4+1)^2}\ge\frac{1}{2}-\frac{1+2}{2(1+4)}=\frac{1}{5}.
\end{multline}

Now we prove that if $\alpha$ is a Pisot number of degree not greater than 4 and $\alpha$ is not the root of $x^4-x^3-1=0$, $x^3-x-1=0$ and $x^2-x-1=0$ then

\begin{multline}
L(\alpha)\ge \max \{ \left \| \frac{a_1+a_2+a_3+a_4-2}{2(a_1+a_2+a_3+a_4-1)} \right \|, \left \| \frac{a_1-a_2+a_3-a_4}{2(a_1-a_2+a_3-a_4+1)} \right \|,  
\\ \frac {(a_1-a_3)^2-|a_1-a_3|+(a_2-a_4+1)^2-|a_2-a_4+1|}{2(a_1-a_3)^2+2(a_2-a_4+1)^2}\}\ge \frac{4}{13}. 
\end{multline}

It's easy to see that $L(\alpha)=\frac{1}{5}$ if and only if $(a_1, a_2, a_3, a_4)$ satisfy one of the following systems of equations

\begin{multicols}{2}
\begin{equation}
\begin{cases}
a_1+a_2+a_3+a_4=0, \\
a_1+a_3=a_2+a_4, \\
|a_1-a_3|=1, \\
|a_2-a_4+1|=2.
\end{cases}
\end{equation}
\begin{equation}
\begin{cases}
a_1+a_2+a_3+a_4=0, \\
a_1+a_3=a_2+a_4, \\
|a_1-a_3|=2, \\
|a_2-a_4+1|=1.
\end{cases}
\end{equation}
\begin{equation}
\begin{cases}
a_1+a_2+a_3+a_4=2, \\
a_1+a_3=a_2+a_4, \\
|a_1-a_3|=1, \\
|a_2-a_4+1|=2.
\end{cases}
\end{equation}
\begin{equation}
\begin{cases}
a_1+a_2+a_3+a_4=2, \\
a_1+a_3=a_2+a_4, \\
|a_1-a_3|=2, \\
|a_2-a_4+1|=1.
\end{cases}
\end{equation}

\begin{equation}
\begin{cases}
a_1+a_2+a_3+a_4=0, \\
a_1+a_3-a_2-a_4=-2, \\
|a_1-a_3|=1, \\
|a_2-a_4+1|=2.
\end{cases}
\end{equation}
\begin{equation}
\begin{cases}
a_1+a_2+a_3+a_4=0, \\
a_1+a_3-a_2-a_4=-2, \\
|a_1-a_3|=2, \\
|a_2-a_4+1|=1.
\end{cases}
\end{equation}
\begin{equation}
\begin{cases}
a_1+a_2+a_3+a_4=2, \\
a_1+a_3-a_2-a_4=-2, \\
|a_1-a_3|=1, \\
|a_2-a_4+1|=2.
\end{cases}
\end{equation}
\begin{equation}
\begin{cases}
a_1+a_2+a_3+a_4=2, \\
a_1+a_3-a_2-a_4=-2, \\
|a_1-a_3|=2, \\
|a_2-a_4+1|=1.
\end{cases}
\end{equation}
\end{multicols}

System (17) has 4 solutions $(\frac12,\frac12,-\frac12,-\frac12)$, $(\frac12,-\frac32,-\frac12,\frac32)$, $(-\frac12,\frac12,\frac12,-\frac12)$, $(-\frac12,-\frac32,\frac12,\frac32)$.

System (18) has 4 solutions $(1,0,-1,0)$, $(1,-1,-1,1)$, $(-1,0,1,0)$, $(-1,-1,1,1)$.

System (19) has 4 solutions $(1,1,0,0)$, $(1,-1,0,2)$, $(0,1,1,0)$, $(0,-1,1,2)$.

System (20) has 4 solutions $(\frac32,\frac12,-\frac12,\frac12)$, $(\frac32,-\frac12,-\frac12,\frac32)$, $(-\frac12,\frac12,\frac32,\frac12)$, $(-\frac12,-\frac12,\frac32,\frac32)$.

System (21) has 4 solutions $(0,1,-1,0)$, $(0,-1,-1,2)$, $(-1,1,0,0)$, $(-1,-1,0,2)$.

System (22) has 4 solutions $(\frac12,\frac12,-\frac32,\frac12)$, $(\frac12,-\frac12,-\frac32,\frac32)$, $(-\frac32,\frac12,\frac12,\frac12)$, $(-\frac32,-\frac12,\frac12,\frac32)$.

System (23) has 4 solutions $(\frac12,\frac32,-\frac12,\frac12)$, $(\frac12,-\frac12,-\frac12,\frac52)$, $(-\frac12,\frac32,\frac12,\frac12)$, $(-\frac12,-\frac12,\frac12,\frac52)$.

System (24) has 4 solutions $(1,1,-1,1)$, $(1,0,-1,2)$, $(-1,1,1,1)$, $(-1,0,1,2)$.

Among these solutions only two correspond to Pisot numbers, namely $(1,1,0,0)$ and $(0,1,1,0)$.

If $(a_1, a_2, a_3, a_4)$ doesn't satisfy any of the systems above then $(A,B)$ with minimal $A^2+B^2$ are $(2,3)$, $(3,2)$. Then for $z$, defined in (15), we have

$$z\ge\frac{1}{2}-\frac{2+3}{2(4+9)}=\frac{4}{13}.$$

Now we search for Pisot numbers of degree not greater than 4 with $L(\alpha)=\frac{4}{13}$. It happens if and only if $(a_1, a_2, a_3, a_4)$ satisfy one of the following systems of equations

\begin{multicols}{2}
\begin{equation}
\begin{cases}
a_1+a_2+a_3+a_4=0, \\
a_1+a_3=a_2+a_4, \\
|a_1-a_3|=2, \\
|a_2-a_4+1|=3.
\end{cases}
\end{equation}
\begin{equation}
\begin{cases}
a_1+a_2+a_3+a_4=0, \\
a_1+a_3=a_2+a_4, \\
|a_1-a_3|=3, \\
|a_2-a_4+1|=2.
\end{cases}
\end{equation}
\begin{equation}
\begin{cases}
a_1+a_2+a_3+a_4=2, \\
a_1+a_3=a_2+a_4, \\
|a_1-a_3|=2, \\
|a_2-a_4+1|=3.
\end{cases}
\end{equation}
\begin{equation}
\begin{cases}
a_1+a_2+a_3+a_4=2, \\
a_1+a_3=a_2+a_4, \\
|a_1-a_3|=3, \\
|a_2-a_4+1|=2.
\end{cases}
\end{equation}

\begin{equation}
\begin{cases}
a_1+a_2+a_3+a_4=0, \\
a_1+a_3-a_2-a_4=-2, \\
|a_1-a_3|=2, \\
|a_2-a_4+1|=3.
\end{cases}
\end{equation}
\begin{equation}
\begin{cases}
a_1+a_2+a_3+a_4=0, \\
a_1+a_3-a_2-a_4=-2, \\
|a_1-a_3|=3, \\
|a_2-a_4+1|=2.
\end{cases}
\end{equation}
\begin{equation}
\begin{cases}
a_1+a_2+a_3+a_4=2, \\
a_1+a_3-a_2-a_4=-2, \\
|a_1-a_3|=2, \\
|a_2-a_4+1|=3.
\end{cases}
\end{equation}
\begin{equation}
\begin{cases}
a_1+a_2+a_3+a_4=2, \\
a_1+a_3-a_2-a_4=-2, \\
|a_1-a_3|=3, \\
|a_2-a_4+1|=2.
\end{cases}
\end{equation}

\end{multicols}

One can easily see that the only solution of these systems that corresponds to Pisot number is $(2,1,-1,0)$. So if $z\ne 0,\frac{1}{5}, \frac{4}{13}$ then

$$\frac {(a_1-a_3)^2-|a_1-a_3|+(a_2-a_4+1)^2-|a_2-a_4+1|}{2(a_1-a_3)^2+2(a_2-a_4+1)^2}=$$
$$=\frac{1}{2}-\frac{|a_1-a_3|+|a_2-a_4+1|}{2(a_1-a_3)^2+2(a_2-a_4+1)^2}\ge\frac{1}{2}-\frac{1+4}{2(1+16)}=\frac{6}{17}>\frac{1}{3}.$$

Theorem 4 is proved.

\section{Pisot numbers less than $\frac{\sqrt 5 + 1}{2}$}

Here we prove Theorem 5.

From Proposition 1.2 it follows that it's enough to consider several cases. 

Let $p\ge1$, $p \equiv 0 \pmod 2$ and $\alpha$ be the maximal zero of $T_{2p+1}=\frac{x^{2p+1}(x^2-x-1)+1}{x^2-1}$. It's easy to see that $T_{2p+1}=x^{2p+1}-\sum_{i=0}^{p} x^{2i}$. Then from Theorem B we have
$L(\alpha)=\frac{1}{2}.$

Let $p\ge1$, $p \equiv 1 \pmod 2$ and $\alpha$ be the maximal zero of $T_{2p+1}=\frac{x^{2p+1}(x^2-x-1)+1}{x^2-1}$. Then from Theorem C we have
$L(\alpha)=\frac{p+1}{2p+4}\ge \frac{1}{3}.$

Let $p\ge2$ and $\alpha$ be the maximal zero of $T_{2p}=\frac{x^{2p}(x^2-x-1)+1}{x-1}$. It's easy to see that $T_{2p}=x^{2p+1}-\sum_{i=0}^{2p-1}x^i$. Then from Theorem 2 we have $L(\theta_{\alpha})\ge\frac{p-1}{2p-1}\ge\frac{1}{3}.$

Let $\alpha$ be the maximal root of $x^{k+2}-x^{k+1}-x^{k}+x^2-1=0$.

If $k=1$, then $L(\alpha)=\frac{1}{5}$. If $k=2$, then $L(\alpha)=\frac{3}{17}$.

If $k \equiv 0,1,3 \pmod 4$ we deal with periodic sequence of length 4 from Lemma 3.4. If $k \equiv 2 \pmod 4$ the sequence from Lemma 3.4 has a zero element. So it gives no non-zero lower bound for $L(\alpha)$. That's why in the case $k \equiv 2 \pmod 4$ we construct periodic sequences of length 8. We construct these sequences just for the special recurrence associated with the polynomial $x^{k+2}-x^{k+1}-x^{k}+x^2-1=0$.

If $k \equiv 3\pmod 4$, then $s_{1,4}=2$, $s_{2,4}=1$, $s_{3,4}=-1$, $s_{4,4}=0$. It's easy to see that the period of the sequence from Lemma 3.4 is $$\frac{4}{13}, \frac{6}{13}, \frac{9}{13}, \frac{7}{13}.$$ Therefore, $L(\alpha)\ge \frac{4}{13}$.

If $k \equiv 0\pmod 4$, then $s_{1,4}=1$, $s_{2,4}=2$, $s_{3,4}=0$, $s_{4,4}=-1$. It's easy to see that the period of the sequence from Lemma 3.4 is $$\frac{6}{17}, \frac{10}{17}, \frac{11}{17}, \frac{7}{17}.$$ Therefore, $L(\alpha)\ge \frac{6}{17}$.

If $k \equiv 1\pmod 4$, then $s_{1,4}=0$, $s_{2,4}=1$, $s_{3,4}=1$, $s_{4,4}=0$. It's easy to see that the period of the sequence from Lemma 3.4 is $$\frac{1}{5}, \frac{2}{5}, \frac{4}{5}, \frac{3}{5}.$$ Therefore, $L(\alpha)\ge \frac{1}{5}$.

If $k \equiv 2\pmod 8$, then $s_{1,8}=1$, $s_{2,8}=0$, $s_{3,8}=0$, $s_{4,8}=1$, $s_{5,8}=0$, $s_{6,8}=0$, $s_{7,8}=0$, $s_{8,8}=0$. According to Lemma 3.1 and Remark from Section 3 for all $n \equiv 2\pmod 8$ there exist periodic sequence modulo 1 with period of length 8 $$\frac{3}{17},\frac{10}{17},\frac{5}{17},\frac{11}{17},\frac{14}{17},\frac{7}{17},\frac{12}{17},\frac{6}{17}.$$ Therefore, $L(\alpha)\ge \frac{3}{17}$.

If $k \equiv 6\pmod 8$, then $s_{1,8}=1$, $s_{2,8}=1$, $s_{3,8}=0$, $s_{4,8}=0$, $s_{5,8}=0$, $s_{6,8}=-1$, $s_{7,8}=0$, $s_{8,8}=1$. According to Lemma 3.1 and Remark from Section 3 for all $n \equiv 6\pmod 8$ there exist periodic sequence modulo 1 with period of length 8 $$\frac{3}{17},\frac{11}{17},\frac{12}{17},\frac{10}{17},\frac{14}{17},\frac{6}{17},\frac{5}{17},\frac{7}{17}.$$ Therefore, $L(\alpha)\ge \frac{3}{17}$.

If $\alpha$ is the maximal root of $1-x+x^2-x^4+2x^5-x^6=0$. Then there exist periodic sequence modulo 1 with period of length 16 $$\frac{4}{17},\frac{7}{17},\frac{11}{17},\frac{12}{17},\frac{12}{17},\frac{8}{17},\frac{12}{17},\frac{7}{17},\frac{13}{17},\frac{10}{17},\frac{6}{17},\frac{5}{17},\frac{5}{17},\frac{9}{17},\frac{5}{17},\frac{10}{17}.$$ Therefore, $L(\alpha)\ge \frac{4}{17}$.

Theorem 5 is proved.

\section{Pisot numbers of degree 2 and 3}

\subsection{Set of coefficients of Pisot numbers of degree 1 and 2}

Consider the sets  \begin{multline} \notag
\Lambda_1=\{a_1 \in \mathbb Z: P(x)= x-a_1 \\ \text{is the minimal polynomial of a Pisot number of degree 1}\},
\end{multline}
\begin{multline} \notag
\Lambda_2=\{(a_1, a_2) \in \mathbb Z^2: P(x)= x^2-a_1x-a_2 \\ \text{is the minimal polynomial of a Pisot number of degree 2}\}.
\end{multline}

The following statement holds.
\\

\noindent {\bf Proposition 6.1} {\it We have
$$\Lambda_1=\{a_1 \in \mathbb Z: a_1>1 \},$$
$$\Lambda_2=\{(a_1, a_2) \in \mathbb Z^2: a_2>0, a_2<a_1+1 \} \cup \{(a_1, a_2) \in \mathbb Z^2: a_2<0, a_2>-a_1+1 \}.$$

}

The set $\Lambda_2$ is shown in Fig. 1. Notice that the ray $a_2=0$, $a_1>1$ corresponds to the set of Pisot numbers of degree 1. 

\begin{figure}[htbp]
\centering
\includegraphics[width=100mm]{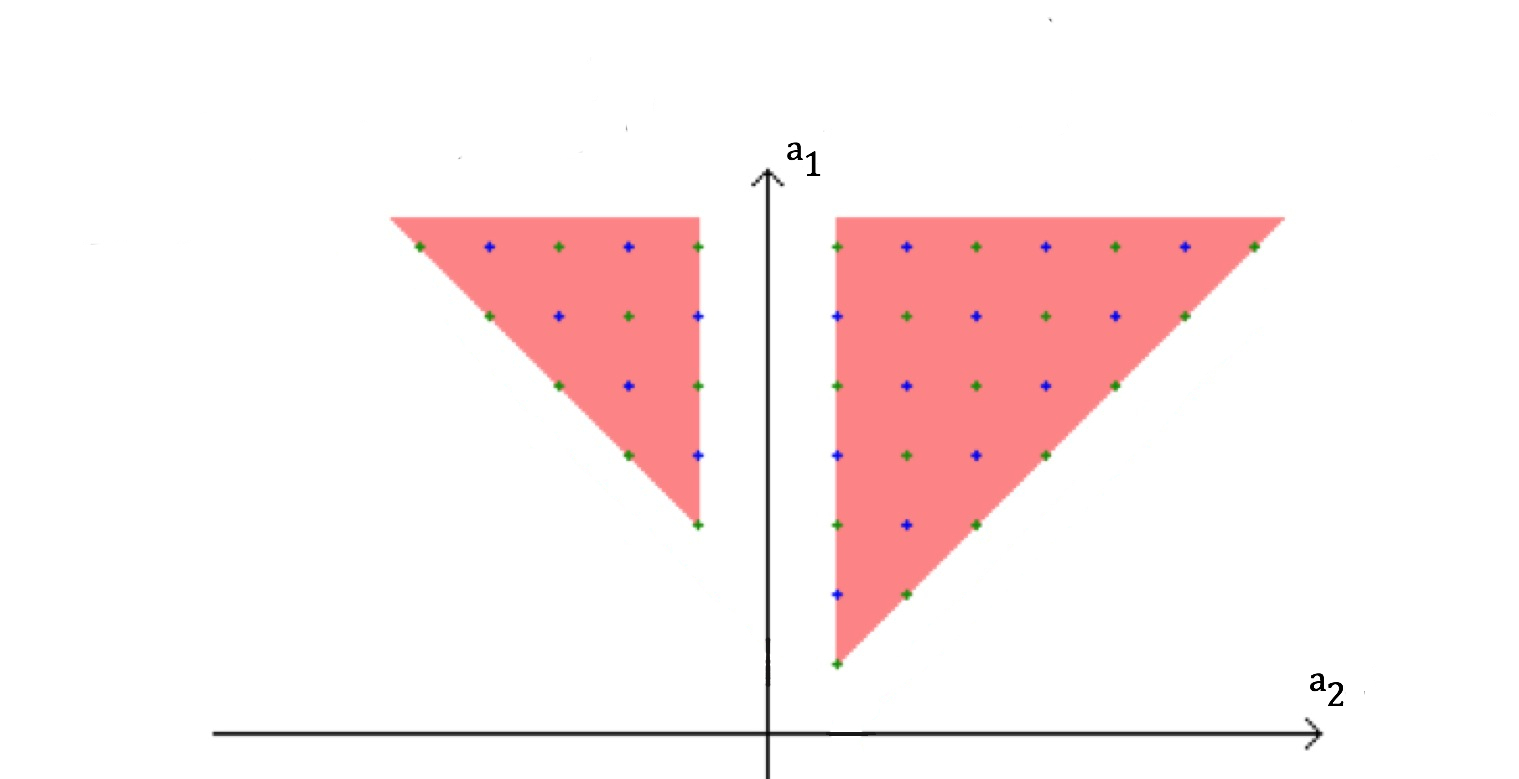}
\caption{Pisot numbers of degree 2}
\label{labelname}
\end{figure}

Due to Theorems A, C and D we know the exact value of $L(\alpha)$ for all Pisot numbers of degree 1 and 2.
\\

\noindent {\bf Proposition 6.2} {\it

1) If $\alpha$ is a Pisot number of degree 1, then $$L(\alpha)=\frac{a_1}{2a_1+2};$$

2) If $\alpha=\frac{\sqrt 5 +1}{2}$, then $$L(\alpha)=\frac{1}{5};$$

3) If $\alpha$ is a Pisot number of degree 2, $\alpha\ne\frac{\sqrt 5 +1}{2}$ and $a_2>0$, then $$L(\alpha)=\frac{a_1+a_2-2}{2a_1+2a_2-2};$$

4) If $\alpha$ is a Pisot number of degree 1 and $a_2<0$, then $$L(\alpha)=\frac{a_1-a_2}{2a_1-2a_2+2}.$$

}

In fact, Proposition 6.2 follows from results of the papers \cite{ref12} by Dubickas, \cite{ref13} and \cite{ref15} by Zhuravleva. 

\subsection{Set of coefficients of Pisot numbers of degree 3}

Let $P(x)=x^3-a_1x^2-a_2x-a_3=0$ be a polynomial with zeros $\alpha_1=\alpha$, $\alpha_2$, $\alpha_3$. Then  $x_n=a_1x_{n-1}+a_2x_{n-2}+a_3x_{n-3}$ and
\begin{equation}
\begin{cases}
a_1=\alpha_1+\alpha_2+\alpha_3, \\
a_2=-\alpha_1 (\alpha_2+\alpha_3)-\alpha_2\alpha_3, \\
a_3=\alpha_1\alpha_2\alpha_3.
\end{cases}
\end{equation}

We consider the set
\begin{multline} \notag
\Lambda_3(a_1)=\{(a_2, a_3) \in \mathbb Z^2: P(x)= x^3-a_1x^2-a_2x-a_3 \\ \text{is the minimal polynomial of a Pisot number of degree 3}\}.
\end{multline}

We rewrite (33) in the following form.
\begin{equation}
\begin{cases}
\alpha_1=a_1-(\alpha_2+\alpha_3), \\
a_2=-a_1(\alpha_2+\alpha_3)+(\alpha_2+\alpha_3)^2 -\alpha_2\alpha_3,\\
a_3=(a_1-(\alpha_2+\alpha_3))\alpha_2\alpha_3.
\end{cases}
\end{equation}

From the definition of Pisot number we have $\alpha_1>1$, $|\alpha_2|<1$, $|\alpha_3|<1$. So (34) leads to
\begin{equation}
\begin{cases}
a_1\ge 0, \\
|a_2|<2a_1+5, \\
|a_3|<a_1+2.
\end{cases}
\end{equation}

We note that if $a_3=0$, then $P(X)$ is reducible. Therefore, we don't consider the pairs $(a_2, 0)$. 

For $0\le a_1 \le 2$ we find all elements of the set $\Lambda_1(a_1)$ using computer. It turned out that these elements are
$(a_1, a_2, a_3)=(0,1,1)$,$(1,0,1)$, $(1,1,1)$, $(1,2,1)$, $(1,2,2)$, $(1,3,2)$, $(2, 1,-1)$, $(2,-1, 1)$, $(2, 0, 1)$, $(2,1, 1)$,$(2,2, 1)$, $(2,3, 1)$, $(2, 0, 2)$, $(2, 1, 2)$, $(2, 2, 2)$, $(2, 3, 2)$, $(2, 4, 2)$, $(2, 3, 3)$, $(2, 4, 2)$, $(2, 5, 2)$.
\\

{\bf Proposition 6.3}
{\it We have} $$\Lambda_3(a_1)=\{(a_2, a_3) \in \mathbb Z^2: a_1+a_2+a_3> 1, a_1-a_2+a_3> -1, a_2> a_3^2-a_1a_3-1, a_3\ne 0\}.$$

Optimal bounds for the coefficients $(a_1, a_2, a_3)$ follow from equalities (34) and $\alpha_1>1$, $|\alpha_2|<1$, $|\alpha_3|<1$ by easy calculations. To complete the proof of Proposition 6.3 it remains to verify the following Lemma. 
 \\

{\bf Lemma 6.1} {\it If $(a_2, a_3) \in \Lambda_3(a_1)$, then $P(X)=x^3-a_1x^2-a_2x-a_3$ is irreducible.}

\begin{proof}

We assume that $P(X)=(x-m_1)(x-m_2)(x-m_3)$ or $P(X)=(x-m_1)(x^2-m_2x-m_3)$, where $m_1, m_2, m_3 \in \mathbb Z$ by Gauss' lemma.

Consider the first case. Without loss of generality we suppose that $m_1>1$, $|m_2|, |m_3|<1$. Then $m_2=m_3=0$. Therefore, $a_3=a_2=0$.

Consider the second case. Let $\beta_1$, $\beta_2$ be the zeros of $x^2-m_2x-m_3$. 

If $|\beta_1|<1$ and $|\beta_2|<1$, then we have $|m_3|=|\beta_1\beta_2|<1$. Therefore, $m_3=0$. It's easy to see that then $m_2=\beta_1=\beta_2=a_2=a_3=0$. 

If $|\beta_1|>1$ and $|\beta_2|<1$, then $|m_1|<1$. Therefore, $m_1=a_3=0$. 

\end{proof}
The figure below represents the set $\Lambda_3(a_1)$ for $a_1=9$. 

Point $A$ corresponds to $\alpha_2=1, \alpha_3=1$.

Point $B$ corresponds to $\alpha_2=1, \alpha_3=-1$.

Point $C$ corresponds to $\alpha_2=-1, \alpha_3=-1$.

Point $E$ corresponds to $\alpha_2=1, \alpha_3=0$.

Point $F$ corresponds to $\alpha_2=0, \alpha_3=-1$.

Point $J$ corresponds to $\alpha_2=i, \alpha_3=-i$.

Point $H$ corresponds to  $\alpha_2=0,5, \alpha_3=0,5$.

Point $I$ corresponds to  $\alpha_2=-0,5, \alpha_3=-0,5$.

Points $A, E, B$ lie on the straight line $a_1+a_2+a_3=1$.

Points $C, F, B$ lie on the straight line $a_1-a_2+a_3= -1$.

Points $A, J, C$ lie on the curve $a_2= a_3^2-a_1a_3-1$.

$\Sigma$ is the curve of zero discriminant. It separates polynomials with three real roots and polynomials with one real root and two complex conjugate roots.

We note that the straight line $a_3=0$ corresponds to Pisot numbers of degree 1 and 2.

\begin{figure}[htbp]
\centering
\includegraphics[width=140mm]{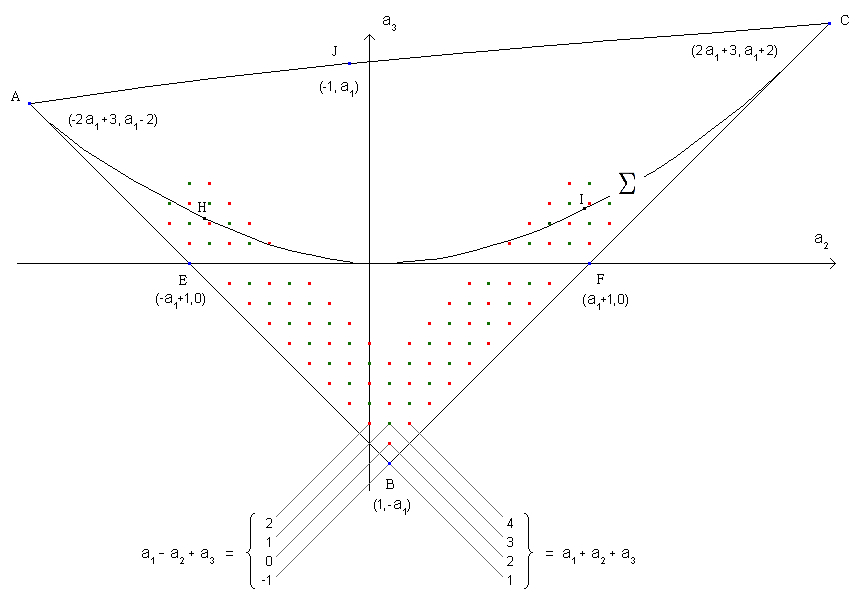}
\caption{Pisot numbers of degree 3}
\label{labelname}
\end{figure}

\newpage

\subsection{Partition of $\Lambda_3(a_1)$ into three domains $\Gamma_1(a_1)$, $\Gamma_2(a_1)$, $\Gamma_3(a_1)$}

Due to Theorem 2 and Proposition 6.3 we divide $\Lambda_3(a_1)$ into three sets according to which function from Theorem 2 is maximal. 

On the figure below

1) the domain $\Gamma_1(a_1)$, where maximum in (2) attains at $\frac{a_1-a_2+a_3}{2(a_1-a_2+a_3+1)}$, is marked in blue;

2) the domain $\Gamma_2(a_1)$, where maximum in (2) attains at $\frac{a_1+a_2+a_3-2}{2(a_1+a_2+a_3-1)}$, is marked in red;

3) the domain $\Gamma_3(a_1)$, where maximum in (2) attains at $\frac {(a_1-a_3)^2-|a_1-a_3|+(a_2+1)^2-|a_2+1|}{2(a_1-a_3)^2+2(a_2+1)^2}$, is marked in green.

\begin{figure}[htbp]
\centering
\includegraphics[width=170mm]{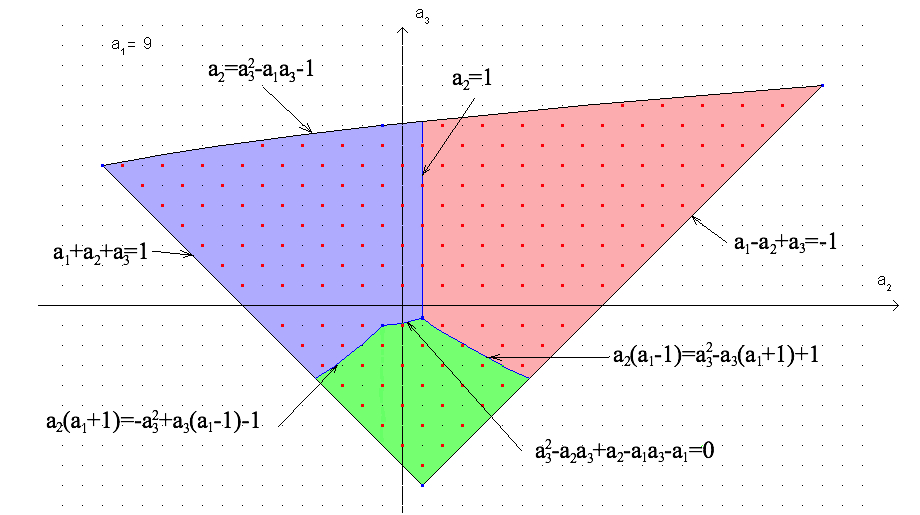}
\caption{Lower bounds for $L(\alpha)$}
\label{labelname}
\end{figure}

We note that in the case $a_3>0$ due to Theorems C and D we have equalities for $L(\alpha)$, namely
\\

\noindent{\bf Corollary 6.1} {\it Let $\alpha$ be a Pisot number of degree 3. Then

1) If $a_3\ge0, a_2\le0$, then $$L(\alpha)=\frac{a_1-a_2+a_3}{2(a_1-a_2+a_3+1)};$$

2) If $a_3>0, a_2>0$, then} $$L(\alpha)=\frac{a_1+a_2+a_3-2}{2(a_1+a_2+a_3-1)}.$$

Theorem 6 gives equalities for the values of $L(\alpha)$ for certain parts of blue and red regions in the case $a_3<0$. For the rest Pisot numbers of degree 3 we don't have exact values for $L(\alpha)$, but a lower bound only.

\subsection{Upper bounds}

Here we prove Theorem 6.

It is enough to verify two following statements.

1) If $a_2<0, a_3<0$ and $a_2+\frac {a_3^2-a_3(a_1-1)+1}{a_1+1}\le 0$, then $$L(\alpha) \le \frac {a_1-a_2+a_3}{2(a_1-a_2+a_3+1)}.$$ 

2) If $a_2>0, a_3<0$ and $a_2-\frac {a_3^2-a_3(a_1+1)+1}{a_1-1}\ge 0$, then
$$ L(\alpha) \le \frac {a_1+a_2+a_3-2}{2(a_1+a_2+a_3-1)}.$$.

Put $u_n:=[\xi x_n ]$ and $w_n := \{\xi x_n \} - \frac{1}{2}$. Then $-\frac{1}{2} \le w_n < \frac{1}{2}$. Therefore, 
$$u_n+w_n+\frac{1}{2}=\sum_{i=1}^{3} a_i(u_{n-i}+w_{n-i}+\frac{1}{2}).$$

\noindent{\bf Lemma 7.1}
{\it We have the following inequalities

1) In case 1 
\begin{equation}
a_2-a_3\le0,
\end{equation}

2) In case 2
\begin{equation}
-a_2-a_3\le-1.
\end{equation}}

\begin{proof}

It is easy to see that $a_2=a_3$ is the tangent to the curve $(a_1+1)a_2+a_3^2-a_3(a_1-1)+1= 0$ at $a_2=-1$, $a_3=-1$. That's why there is no intersection of our curve and the area where $a_2-a_3>0$.

For the curve $(a_1-1)a_2-a_3^2+a_3(a_1-1)-1= 0$ we have the tangent $a_2=-a_3$ at $a_2=-1$, $a_3=1$. But this point doesn't belong to the area we consider in case 2. That's why we have $-a_2-a_3\le -1$.
\end{proof}

We assume that $L(\alpha) > \frac {a_1-a_2+a_3}{2(a_1-a_2+a_3+1)}$ (in case 1) and $L(\alpha)>\frac {a_1+a_2+a_3-2}{2(a_1+a_2+a_3-1)}$ (in case 2). Then there exists a positive number
$$
w<\frac{1}{2}-L(\alpha)=\frac{1}{2}-\frac {a_1-a_2+a_3}{2(a_1-a_2+a_3+1)}=\frac {1}{2(a_1-a_2+a_3+1)} \quad \text{(in case 1)},
$$
$$
w<\frac{1}{2}-L(\alpha)=\frac{1}{2}-\frac {a_1+a_2+a_3-2}{2(a_1+a_2+a_3-1)}=\frac {1}{2(a_1+a_2+a_3-1)}  \quad \text{(in case 2)},
$$
such that $-w \le w_n \le w$ for $n \ge n_0$.
\\

\noindent {\bf Lemma 7.2.}
{\it We have $\delta_{n}:= \sum_{i=1}^{d} a_i w_{n-i}-w_n \in \{-\frac{1}{2}, \frac{1}{2} \}$ for $n\ge n_0$.}

\begin{proof}

We should note that $\delta_n$ is of the form $\frac{s}{2}$, where $s$ is an odd number. Assume that $|\delta_n|\ge\frac{3}{2}$. Then
$$ \frac{3}{2}\le |a_1 w_{n-1}+a_2 w_{n-2}+a_3 w_{n-3}-w_n|\le(|a_1|+|a_2|+|a_3|+1)w.$$

In the first case $$ \frac{3}{2} < \frac {(a_1-a_2-a_3+1)}{2(a_1-a_2+a_3+1)}.$$

Then $a_1-a_2+2a_3+1<0.$ This contradicts the inequalities $a_1+a_2+a_3>1$, $a_1-a_2+a_3>-1$ and (36).

In the second case $$ \frac{3}{2} < \frac {(a_1+a_2-a_3+1)}{2(a_1+a_2+a_3-1)}.$$

Then $a_1+a_2+2a_3-2<0.$ This contradicts the inequalities $a_1+a_2+a_3>1$, $a_1-a_2+a_3>-1$ and (37).

\end{proof}

We take large number $N$ and $n>n_0 + 3$. Put $l_j=1$, if $\delta_{n+j}=\frac{1}{2}$, and $l_j=-1$, if $\delta_{n+j}=-\frac{1}{2}$, where $j=0, \ldots, N-1$. We consider the sum
$$
\sum_{j=0}^{N-1} l_j (\sum_{i=1}^{3} a_i w_{n-i+j}-w_{n+j})=\frac{N}{2}.
$$

If we put $l_j:=0$ for $j<0$ and $j>N-1$, we have
$$
\sum_{j=0}^{N-1} l_j (\sum_{i=1}^{3} a_i w_{n-i+j}-w_{n+j})=\sum_{j=-3}^{N-1} w_{n+j} (\sum_{i=1}^{3}a_i l_{j+i}-l_j).
$$

We estimate $$T:=\sum_{j=-3}^{N-1}|\sum_{i=1}^{3}a_i l_{j+i}-l_j| \le \sum_{j=0}^{N-4}|\sum_{i=1}^{3}a_i l_{j+i}-l_j|+6(\sum_{i=1}^{3}|a_i| + 1).$$

Now the proof for the second case splits into three following cases.

1) For $a_2<0$ we prove that $L(\alpha)=\frac{a_1-a_2+a_3}{2(a_1-a_2+a_3+1)}$;

2) For $a_2>0$, $a_1-a_2+a_3>0$ we prove that $L(\alpha)=\frac{a_1+a_2+a_3-2}{2(a_1+a_2+a_3-1)}$;

2*) For $a_2>0$, $a_1-a_2+a_3=0$ we prove that $L(\alpha)=\frac{a_1+a_2+a_3-2}{2(a_1+a_2+a_3-1)}$.
\\

Consider all vectors $(l_j, l_{j+1}, l_{j+2}, l_{j+3})$ for $N-4\ge j\ge0$. There are 16 possibilities for the vector $(l_j, l_{j+1}, l_{j+2}, l_{j+3})$ to attain. All possible cases are listed in the table below, where + stands for $1$, and - stands for $-1$. Put $R_j=|a_1 l_{j+1}+a_2 l_{j+2}+a_3 l_{j+3}-l_j|$, $U=a_1-a_2+a_3+1$ and $V=a_1+a_2+a_3-1$. We use the values marked by * only for case 2*.

\begin{center}
Table 1
\end{center}
\noindent
\begin{tabular}{|l|l|l|l|l|l|l|l|l|l|}
\hline
Type & $l_j$ & $l_{j+1}$ & $l_{j+2}$ & $l_{j+3}$ & $R_j$ & $R_j-U$ & & $R_j-V$ & \\
\hline
A1 & - & + & + & + & $a_1+a_2+a_3+1$ & $2a_2$ & $\le 0$ & 2 & $>0$ \\
\hline
A2 & + & - & - & - & $a_1+a_2+a_3+1$ & $2a_2$ & $\le 0$ & 2 & $>0$\\
\hline
B1 & + & + & + & + & $a_1+a_2+a_3-1$ & $2a_2-2$ & $\le 0$ & 0 & $\le0$ \\
\hline
B2 & - & - & - & - & $a_1+a_2+a_3-1$ & $2a_2-2$ & $\le 0$ & 0 & $\le0$\\
\hline
C1 & - & + & - & - & $a_1-a_2-a_3+1$ & $-2a_3$ & $> 0$ & $-2a_2-2a_3+2$ & $\le0$\\
\hline
C2 & + & - & + & + & $a_1-a_2-a_3+1$ & $-2a_3$ & $> 0$ & $-2a_2-2a_3+2$ & $\le0$\\
\hline
D1 & + & + & - & - & $a_1-a_2-a_3-1$ & $-2a_3-2$ & $> 0$ & $-2a_2-2a_3$ & $\le0$\\
\hline
D2 & - & - & + & + & $a_1-a_2-a_3-1$ & $-2a_3-2$ & $> 0$ & $-2a_2-2a_3$ & $\le0$\\
\hline
E1 & - & + & + & - & $a_1+a_2-a_3+1$ & $2a_2-2a_3$ & $\le 0$ & $-2a_3+2$ & $>0$\\
\hline
E2 & + & - & - & + & $a_1+a_2-a_3+1$ & $2a_2-2a_3$ & $\le 0$ & $-2a_3+2$ & $>0$\\
\hline
F1 & + & + & + & - & $a_1+a_2-a_3-1$ & $2a_2-2a_3-2$ & $\le 0$ & $-2a_3$ & $>0$\\
\hline
F2 & - & - & - & + & $a_1+a_2-a_3-1$ & $2a_2-2a_3-2$ & $\le 0$ & $-2a_3$ & $>0$\\
\hline
G1 & - & + & - & + & $a_1-a_2+a_3+1$ & $0$ & $\le 0$ & $-2a_2+2$ & $\le0$\\
\hline
G2 & + & - & + & - & $a_1-a_2+a_3+1$ & $0$ & $\le 0$ & $-2a_2+2$ & $\le0$\\
\hline
H1 & + & + & - & + & $a_1-a_2+a_3-1$ or $1$*  & $-2$ & $\le 0$ & $-2a_2$ or $-2a_2+2$* & $\le0$\\
\hline
H2 & - & - & + & - & $a_1-a_2+a_3-1$ or $1$* & $-2$ & $\le 0$ & $-2a_2$ or $-2a_2+2$* & $\le0$\\
\hline
\end{tabular}
\\

Let the sum T contain $|A|$ elements of types $A_1$ and $A_2$, $|B|$ elements of types $B_1$ and $B_2$, $|C|$ elements of types $C_1$ and $C_2$, $|D|$ elements of types $D_1$ and $D_2$, $|E|$ elements of types $E_1$ and $E_2$, $|F|$ elements of types $F_1$ and $F_2$, $|G|$ elements of types $G_1$ and $G_2$, $|H|$ elements of types $H_1$ and $H_2$. Then
\begin{multline}
T\le |1+a_1+a_2+a_3| |A|+|-1+a_1+a_2+a_3||B|+|1+a_1-a_2-a_3||C|+|-1+a_1-a_2-a_3||D|+\\
+|1+a_1+a_2-a_3||E|+|-1+a_1+a_2-a_3||F|+|1+a_1-a_2+a_3||G|+|-1+a_1-a_2+a_3||H|+
c_1(a_1, a_2, a_3).
\end{multline}

Hereinafter we denote by $c_i(a_1, a_2, a_3)$ some constants that depend on $a_1, a_2, a_3$, but don't depend on $N$.
\\

For our proof we need the following statement.
\\

{\bf Lemma 7.3.} {\it We have}

1) $|C|+|D|\le |E|+|A|+1$;

2) $|E|+|F|\le |H|+|D|+1$;

3) $|A|\le|F|+1$;

4) $|C|\le|H|+1$.

\begin{proof}
1) We note that elements of types $C$ and $D$ are followed by elements of types $E$ and $A$. Therefore, $|C|+|D|\le 
|E|+|A|+1$.

2) Elements of types  $E$ and $F$ are followed by elements of types $H$ and $D$. Therefore, $|E|+|F|\le 
|H|+|D|+1$.

3) This statement was proved in \cite{ref12}.

4) We prove that between any two elements of type $C$ there is an element of type $H$. So, $|C|\le|H|+1$.

We note that an element of type $C$ is preceded by either an element of type $H$, or an element of type $G$. An element of type $G$ is preceded by either an element of type $G$, or an element of type $H$.  For an element of type $C$ we need two + or - standing side by side. This is exactly an element of type $H$.
\end{proof}

\subsubsection{Case $a_2<0$}

We substitute the values from Table 1 into (38) and we get
\begin{multline} \notag
T\le (1+a_1+a_2+a_3) |A|+(-1+a_1+a_2+a_3)|B|+|(1+a_1-a_2-a_3)|C|+(-1+a_1-a_2-a_3)|D|+\\
+(1+a_1+a_2-a_3)|E|+(-1+a_1+a_2-a_3)|F|+(1+a_1-a_2+a_3)|G|+(-1+a_1-a_2+a_3)|H|+
c_1(a_1, a_2, a_3)\le\\
\le(a_1-a_2+a_3+1)N+2a_2|A|+(2a_2-2)|B|+(-2a_3)|C|+(-2a_3-2)|D|+\\
+(2a_2-2a_3)|E|+(2a_2-2a_3-2)|F|+ (-2)|H|+c_2(a_1, a_2, a_3)\le \\ \le (a_1-a_2+a_3+1)N+2a_2|A|+(-2a_3)|C|+(-2a_3-2)|D|+(2a_2-2a_3)|E|+ (-2)|H|+c_2(a_1, a_2, a_3).
\end{multline}

Since from Lemma 7.3 we have $C\le |H|+1$ and $|C|+|D|\le |E|+|A|+1$, then $$(-2a_3-2)(|C|+|D|)+2|C|\le(-2a_3-2)(|E|+|A|+1)+2|H|+2.$$

So,
$$T\le(a_1-a_2+a_3+1)N+(2a_2-2a_3-2)|A|+(2a_2-4a_3-2)|E|+c_3(a_1, a_2, a_3).$$

Dubickas' argument from paper \cite{newdub} enables us to get the following statement.

\noindent {\bf Lemma 7.4.} {\it If $a_2<0, a_3<0$ and $a_2+\frac {a_3^2-a_3(a_1-1)+1}{a_1+1}\le 0$, then $|E|=0$.}

\begin{proof}

Let $(s_1, s_2, s_3, s_4) \in \mathbb Z^4$. From the definition of $l_j$ we have $l_j=2\delta_{n+j}$. We put $v_n=2w_n$. Then 

\begin{multline}
s_1 l_j + s_2 l_{j+1} +s_3 l_{j+2} + s_4 l_{j+3}=\\ = t_1 v_{n+j+3}+t_2 v_{n+j+2}+t_3 v_{n+j+1}+t_4 v_{n+j}+t_5 v_{n+j-1} + t_6 v_{n+j-2}+t_7 v_{n+j-3}.
\end{multline}

Let $$r:=\frac{|s_1 l_j + s_2 l_{j+1} +s_3 l_{j+2} + s_4 l_{j+3}|}{|t_1|+|t_2|+|t_3|+|t_4|+|t_5|+|t_6|+|t_7|}. $$

According to Dubickas \cite{newdub} the pattern $(l_{j}, l_{j+1}, l_{j+2}, l_{j+3})$ doesn't occurs in the sequence $(l_j)_{j=n_0}^{\infty}$ if there exists an integer vector $(s_1, s_2, s_3, s_4) \in \mathbb Z^n$ such that $$r\ge 1-2L(\alpha).$$

This inequality implies that the absolute value of at least one of numbers $v_{n+j-3}$, $v_{n+j-2}$, $v_{n+j-1}$, $v_{n+j}$, $v_{n+j+1}$, $v_{n+j+2}$, $v_{n+j+3}$ must be greater or equal to $r\ge 1-2L(\alpha)$. It's easy to see that
this contradicts $|w_n|\le w$.

Type E deals with patterns $(1,-1,-1,1)$ and $(-1,1,1,-1)$. It's enough to consider the case when $l_j=1$, $l_{j+1}=-1$, $l_{j+2}=-1$, $l_{j+3}=1$. 

Put $s_1=a_3$, $s_2=-a_2$, $s_3=a_1$, $s_4=1$. Then 
$$s_1 l_j + s_2 l_{j+1} +s_3 l_{j+2} + s_4 l_{j+3}=-v_{n+j+3}+(a_1^2+2a_2)v_{n+j+1}+(2a_1a_3-a_2^2)v_{n+j-1}+a_3^2 v_{n+j-3}.$$

Therefore, $$r=\frac{a_1-a_2-a_3-1}{(a_1-a_3)^2+(a_2+1)^2}.$$ 

The inequality $$\frac{a_1-a_2-a_3-1}{(a_1-a_3)^2+(a_2+1)^2}\ge 1-2L(\alpha)$$ is equivalent to the inequality $$\frac{a_1-a_2+a_3}{2(a_1-a_2+a_3+1)} \ge \frac {(a_1-a_3)^2-|a_1-a_3|+(a_2+1)^2-|a_2+1|}{2(a_1-a_3)^2+2(a_2+1)^2}$$ that we use for partition into three domains in subsection 6.3. So in this case this inequality is valid.

\end{proof}

From Lemmas 7.1 and 7.4 we have
$$T\le(a_1-a_2+a_3+1)N+c_3(a_1, a_2, a_3).$$

As $|w_{n+j}| \le w$, then $Tw \ge \frac {N}{2}$. Therefore,
$$
\frac{1}{2w}\le\frac{T}{N} \le \frac{(a_1-a_2+a_3+1)N}{N} +\frac{c_3(a_1, a_2, a_3)}{N}.
$$

This contradicts $w<\frac{1}{2(a_1-a_2+a_3+1)}$, for sufficiently large $N$.

\subsubsection{Case $a_2>0$, $a_1-a_2+a_3>0$}

We substitute the values from Table 1 into (38) and we get
$$T\le(a_1+a_2+a_3-1)N+2|A|+(-2a_2-2a_3+2)|C|+(-2a_2-2a_3)|D|+(-2a_3+2)|E|+(-2a_3)|F|+$$
$$+(-2a_2+2)|G| + (-2a_2)|H|+3(-2a_3+2)+c_4(a_1, a_2, a_3)\le$$
$$\le(a_1+a_2+a_3-1)N+2|A|+(-2a_2-2a_3)|D|+(-2a_3+2)|E|+(-2a_3)|F|+$$
$$ + (-2a_2)|H|+c_5(a_1, a_2, a_3)$$

Since from Lemma 7.3 we have $|E|+|F|\le |D|+|H|+1$ and $|A| \le |F|+1$, then
$$T\le (a_1+a_2+a_3-1)N +(-2a_2-2a_3)|D| + (-2a_2)|H| + (-2a_3+2)(|E|+|F|)+c_5(a_1, a_2, a_3)\le$$
$$\le (a_1+a_2+a_3-1)N+(-2a_2-4a_3+2)|D| + (-2a_2-2a_3+2)|H|+c_6(a_1, a_2, a_3).$$

\noindent {\bf Lemma 7.5.} {\it If $a_2>0, a_3<0$ and $a_2-\frac {a_3^2-a_3(a_1+1)+1}{a_1-1}\ge 0$, then $|D|=0$.}

\begin{proof}

Type D deals with patterns $(1,1,-1,-1)$ and $(-1,-1,1,1)$.It's enough to consider the case when $l_j=1$, $l_{j+1}=1$, $l_{j+2}=-1$, $l_{j+3}=-1$. 

Put $s_1=a_3$, $s_2=-a_2$, $s_3=a_1$, $s_4=1$. Then 
$$s_1 l_j + s_2 l_{j+1} +s_3 l_{j+2} + s_4 l_{j+3}=-v_{n+j+3}+(a_1^2+2a_2)v_{n+j+1}+(2a_1a_3-a_2^2)v_{n+j-1}+a_3^2 v_{n+j-3}.$$

Therefore, $$r=\frac{a_1+a_2-a_3+1}{(a_1-a_3)^2+(a_2+1)^2}.$$ 

The inequality $$\frac{a_1+a_2-a_3+1}{(a_1-a_3)^2+(a_2+1)^2}\ge 1-2L(\alpha)$$ is equivalent to the inequality $$\frac{a_1+a_2+a_3-2}{2(a_1+a_2+a_3-1)} \ge \frac {(a_1-a_3)^2-|a_1-a_3|+(a_2+1)^2-|a_2+1|}{2(a_1-a_3)^2+2(a_2+1)^2}$$ that we use for partition into three domains in subsection 6.3. So in this case this inequality is valid.

\end{proof}

From Lemmas 7.1 and 7.5 we have 
$$T \le (a_1+a_2+a_3-1)N+c_6(a_1, a_2, a_3).$$

As $|w_{n+j}| \le w$, then $Tw \ge \frac {N}{2}$. Therefore,
$$
\frac{1}{2w}\le\frac{T}{N} \le \frac{(\sum_{i=1}^{d} a_1+a_2+a_3-1)N}{N} +\frac{c_6(a_1, a_2, a_3)}{N}.
$$

This contradicts $w<\frac{1}{2(a_1+a_2+a_3-1)}$, for sufficiently large $N$.

\subsubsection{Case $a_2>0$, $a_1-a_2+a_3=0$}

We substitute the values from Table 1 into (38) and we get
$$T\le(a_1+a_2+a_3-1)N+2|A|+(-2a_2-2a_3+2)|C|+(-2a_2-2a_3)|D|+(-2a_3+2)|E|+(-2a_3)|F|+$$
$$+(-2a_2+2)|G| + (-2a_2+2)|H|+3(-2a_3+2)+c_7(a_1, a_2, a_3)\le$$
$$\le (a_1+a_2+a_3-1)N+2|A|+(-2a_2-2a_3)|D|+(-2a_3+2)|E|+(-2a_3)|F|+$$
$$+ (-2a_2+2)|H|+c_8(a_1, a_2, a_3).$$

Since from Lemma 7.3 we have $|E|+|F|\le |D|+|H|+1$ and $|A| \le |F|+1$, then
$$T\le (a_1+a_2+a_3-1)N +(-2a_2-2a_3)|D| + (-2a_2+2)|H| + (-2a_3+2)(|E|+|F|)+c_8(a_1, a_2, a_3)\le$$
$$\le (a_1+a_2+a_3-1)N+(-2a_2-4a_3+2)|D| + (-2a_2-2a_3+4)|H|+c_9(a_1, a_2, a_3).$$

If $a_2+a_3\ge 2$ then using Lemma 7.5 we get $T \le (a_1+a_2+a_3-1)N+c_9(a_1, a_2, a_3).$

If $a_2+a_3= 1$ we need to prove the following result.

\noindent {\bf Lemma 7.6.} {\it If $a_2>0, a_3<0$, $a_2-\frac {a_3^2-a_3(a_1+1)+1}{a_1-1}\ge 0$, $a_1-a_2+a_3=0$ and $a_2+a_3=1$ then $|H|=0$.}

\begin{proof}

It's easy to see that $a_2=\frac{1+a_1}{2}$ and $a_3=\frac{1-a_1}{2}$.

Type H deals with patterns $(1,1,-1,1)$ and $(-1,-1,1,-1)$. It's enough to consider the case when $l_j=1$, $l_{j+1}=1$, $l_{j+2}=-1$, $l_{j+3}=1$. 

Put $s_1=-a_3$, $s_2=a_1$, $s_3=-a_2$, $s_4=a_1-a_3$. Then 
\begin{multline}
s_1 l_j + s_2 l_{j+1} +s_3 l_{j+2} + s_4 l_{j+3}=\\
=\frac{(6a_1-2)v_{n+j+3}+(-6a_1^2-2)v_{n+j+2}+(a_1^2-4a_1-1)v_{n+j+1}}{4}+\\+\frac{(-5a_1^2+1)v_{n+j-1}+(a_1^2-2a_1+1)v_{n+j-2}+(a_1^2-2a_1+1)v_{n+j-3}}{4}.
\end{multline}

We note that if $0<a_1<5$ then there is no $(a_1, a_2, a_3) \in \mathbb Z^3$ that satisfies the conditions of Lemma.

As $a_1\ge 5$ we have $$r=\frac{1}{a}=1-2L(\alpha).$$ 

\end{proof}

From Lemmas 7.1, 7.5 and 7.6 we have
$$T \le (a_1+a_2+a_3-1)N+c_9(a_1, a_2, a_3).$$

As $|w_{n+j}| \le w$, then $Tw \ge \frac {N}{2}$. Therefore,
$$
\frac{1}{2w}\le\frac{T}{N} \le \frac{(\sum_{i=1}^{d} a_1+a_2+a_3-1)N}{N} +\frac{c_9(a_1, a_2, a_3)}{N}.
$$

This contradicts $w<\frac{1}{2(a_1+a_2+a_3-1)}$, for sufficiently large $N$.

\section{Computer calculations}

In some cases it's possible to calculate the exact value of $L(\alpha)$ using computer.

\subsection{Example of calculation from \cite{ref14}}

In paper \cite{ref14} author calculated the value of $L(\alpha)$ for the case $a_1=0, a_2=1, a_3=1$. 
In order to prove the inequality $L(\alpha)\le\frac{1}{5}$ the graph of $F(x,y,z):=\min\{||x||, ||y||, ||z||, ||x+y||, ||y+z||, ||x+y+z||,$$ $$||x+2y+z||, ||x+2y+2z||, ||2x+3y+2z||\}$ for $0\le x,y,z <1$ was carefully studied. It consists of pieces of three-dimensional affine planes. It's obvious that $F(x,y,z)$ can attain its extremums only at the points that are intersections of 4 such planes. 

Using computer calculations, it was found that $F(x,y,z)$ attains its maximal value at 4 following points: $(x_1,y_1,z_1):=(\frac{1}{5}, \frac{2}{5}, \frac{4}{5})$, $(x_2,y_2,z_2):=(\frac{2}{5}, \frac{4}{5}, \frac{3}{5})$, $(x_3,y_3,z_3):=(\frac{4}{5}, \frac{3}{5}, \frac{1}{5})$, $(x_4,y_4,z_4):=(\frac{3}{5}, \frac{1}{5}, \frac{2}{5})$.

\subsection{Results for Pisot numbers of degree 3}

Using the approach described in \cite{ref14} it's possible to calculate the exact value of $L(\alpha)$ for other Pisot numbers of degree 3. Below we list the results of these calculations for all Pisot numbers of degree 3 with $a_3<0$, $2 \le a_1\le 8$ (see Tables 2 - 8) . Each of these 7 tables correspond to one $a_1\in [2,8]$. In the first line twe write all possible values of $a_2$, while in the first column there are all possible values of $a_3$. That's why the intersection of some column and some line corresponds only to one Pisot number $\alpha$. The cells that correspond to such $a_1, a_2, a_3$, that $a_1+a_2+a_3$ is odd, are left empty due to Theorem B which states that $L(\alpha)=\frac{1}{2}$ in this case. In the remaining cells we write the value of $L(\alpha)$, which was obtained with the use of computer calculations.

\begin{center}
Table 2: $a_1=2$

\quad

\begin{tabular}{|l|l|}
\hline
$a_3$ $\backslash$ $a_2$  & 1\\
\hline
-1 & 4/13  \\
\hline
\end{tabular}
\end{center}

\quad

\begin{center}

Table 3: $a_1=2$

\quad

\begin{tabular}{|l|l|l|l|}
\hline
$a_3$ $\backslash$ $a_2$ & 0 & 1 & 2 \\
\hline
-1 & 6/17 & & 9/25 \\
\hline
-2 &  & 11/29 &  \\
\hline
\end{tabular}
\end{center}

\quad

\begin{center}

\newpage

Table 4: $a_1=4$

\quad

\begin{tabular}{|l|l|l|l|l|l|}
\hline
$a_3$ $\backslash$ $a_2$ &-1 & 0 & 1 & 2 & 3\\
\hline
-1 & 2/5 & & 11/29 & & 2/5 \\
\hline
-2 &  & 15/37& &6/15 &  \\
\hline
-3 &  & & 22/53 & &  \\
\hline
\end{tabular}
\end{center}

\begin{center}

Table 5: $a_1=5$

\quad

\begin{tabular}{|l|l|l|l|l|l|l|l|}
\hline
$a_3$ $\backslash$ $a_2$ & -2 &-1 & 0 & 1 & 2 & 3 & 4\\
\hline
-1 & 3/7 & & 15/37 & & 2/5 & & 3/7 \\
\hline
-2 &  & 3/7 &  & 22/53 &  & 27/65 &  \\
\hline
-3 &  & & 28/65 & & 31/73 & &  \\
\hline
-4 &  & &  & 37/85 & & &  \\
\hline
\end{tabular}
\end{center}

\quad

\begin{center}

Table 6: $a_1=6$

\quad

\begin{tabular}{|l|l|l|l|l|l|l|l|l|l|l|l|}
\hline
$a_3$ $\backslash$ $a_2$ &-3& -2 &-1 & 0 & 1 & 2 & 3 & 4 & 5\\
\hline
-1 & 4/9 & & 3/7 & & 22/53 & & 3/7 & & 4/9 \\
\hline
-2 & & 28/65 & & 28/65 &  & 31/73 &  & 3/7 &  \\
\hline
-3 &  & & 4/9 & & 37/85 & & 42/97 & & \\
\hline
-4 &  & & & 45/101 & &48/109 &  & & \\
\hline
-5 &  & & & & 56/125 & &  & &  \\
\hline
\end{tabular}
\end{center}

\quad

\begin{center}

Table 7: $a_1=7$

\quad
\begin{tabular}{|l|l|l|l|l|l|l|l|l|l|l|l|l|l|}
\hline
$a_3$ $\backslash$ $a_2$ &-4 &-3& -2 &-1 & 0 & 1 & 2 & 3 & 4 & 5 & 6\\
\hline
-1 & 5/11 & & 4/9 & & 28/65 & & 3/7 & & 4/9 & & 5/11 \\
\hline
-2 & & 4/9 & & 4/9 &  & 37/85 &  & 42/97 & & 4/9 &   \\
\hline
-3 &  & & 45/101 & & 45/101 & & 11/27 & & 11/25 & & \\
\hline
-4 &  & & & 5/11 & & 56/125 &  & 2/5 & & &\\
\hline
-5 &  & & & & 66/145 & & 23/51 & & & &  \\
\hline
-6 &  & & & &  & 79/173 & & & & &  \\
\hline
\end{tabular}
\end{center}

\quad

\begin{center}

Table 8: $a_1=8$

\quad

\begin{tabular}{|l|l|l|l|l|l|l|l|l|l|l|l|l|l|l|l|}
\hline
$a_3$ $\backslash$ $a_2$ & -5&-4 &-3& -2 &-1 & 0 & 1 & 2 & 3 & 4 & 5 & 6 & 7\\
\hline
-1 & $\frac{6}{13}$ & & $\frac{5}{11}$  & & $\frac{4}{9}$  & & $\frac{37}{85}$  & & $\frac{4}{9}$ & & $\frac{5}{11}$  & & $\frac{6}{13}$ \\
\hline
-2 & & $\frac{5}{11}$  & &$\frac{45}{101}$ &  & $\frac{45}{101}$  &  & $\frac{48}{109}$  & & $\frac{4}{9}$  & & $\frac{5}{11}$  &   \\
\hline
-3 &  & & $\frac{56}{125}$ & & $\frac{5}{11}$  & & $\frac{56}{125}$ & & $\frac{61}{137}$& & $\frac{70}{157}$ & & \\
\hline
-4 &  & & & $\frac{66}{145}$ & &  $\frac{66}{145}$ &  &  $\frac{23}{51}$ & &  $\frac{76}{169}$& & &\\
\hline
-5 &  & & & &  $\frac{6}{13}$ & & $\frac{79}{173}$ & &  $\frac{84}{185}$ & & & &  \\
\hline
-6 &  & & & &  &  $\frac{91}{197}$ & &  $\frac{94}{205}$ & & & & &  \\
\hline
-7 &  & & & &  &  &  $\frac{106}{229}$ & & & & &  &\\
\hline
\end{tabular}
\end{center}

\subsection{Results for Pisot numbers of degree 4}

Using computer it's possible to calculate the exact value of $L(\alpha)$ for some Pisot numbers of degree 4. Calculations are based on the method which is different from the one used in \cite{ref14}. We don't described it in the present paper. Here we just announce some results of our calculations. We plan to publish a detailed description of this algorithm in a separate paper. 

In the table below in the first four columns we write coefficients $a_1, a_2, a_3, a_4$. In the fifth column we put the period of the periodic modulo 1 sequence $X=(x_n)_{n=1}^{\infty}$, which satisfies (1). The length of this period is $t$. For this sequence we have $\min_{i=1, \ldots, t} ||x_i||=L(\alpha)$.The value of $L(\alpha)$ is written in the last column.

\begin{center}

Table 9

\quad

\begin{tabular}{|l|l|l|l|l|l|l|}
\hline
$a_1$ & $a_2$ & $a_3$ & $a_4$ & Sequence & $L(\alpha)$  \\
\hline
2 & 2 & 1 & -1 & $6/17,10/17,11/17,7/17$ & $6/17$ \\
\hline
3 & 0 & 0 & 1 & $22/63, 25/63, 37/63$ & $22/63$ \\
\hline
3 & 3 & -1 & -3 & $27/65, 34/65, 38/65, 31/65$ & $27/65$ \\
\hline
4 & 1 & 1 & -4 & $4/9, 5/9$ & $4/9$ \\
\hline
4 & 8 & 2 & -2 & $5/11$ & $5/11$\\
\hline
5 & 0 & -1 & -2 & $3/7, 3/7, 3/7, 4/7, 4/7, 4/7$ & $3/7$\\
\hline
5 & 1 & -1 & -1 & $87/215, 93/215, 92/215, 128/215, 122/215, 123/215$ & $87/215$\\
\hline
5 & -1 & 0 & 2 & $17/40,17/40,22/40$ & $17/40$ \\
\hline
5 & 0 & 0 & 1 & $87/215,92/215,122/215$ & $87/215$\\
\hline
\end{tabular}
\newline
\end{center}

\end {document}